\RequirePackage{ifpdf}
\ifpdf % We are running pdfTeX in pdf mode
\documentclass[pdftex]{sigma}
\else
\documentclass{sigma}
\fi

\newcommand{\N}{{\mathbb N}}

\newcommand{\R}{{\mathbb R}}

\newcommand{\Z}{{\mathbb Z}}
\newcommand{\C}{{\mathbb C}}

\newcommand{\Rmu}{{\rm Re}\, \mu}
\newcommand{\Rz}{{\rm Re}\, z}

\begin{document}

\allowdisplaybreaks

\renewcommand{\thefootnote}{$\star$}

\renewcommand{\PaperNumber}{050}

\FirstPageHeading

\ShortArticleName{Parameter Dif\/ferentiation for Integral Representations of Associated Legendre Functions}

\ArticleName{On Parameter Dif\/ferentiation for Integral\\ Representations of Associated Legendre Functions\footnote{This paper is a
contribution to the Special Issue ``Symmetry, Separation, Super-integrability and Special Functions~(S$^4$)''. The
full collection is available at
\href{http://www.emis.de/journals/SIGMA/S4.html}{http://www.emis.de/journals/SIGMA/S4.html}}}

\Author{Howard S.~COHL $^{\dag\ddag}$}

\AuthorNameForHeading{H.S.~Cohl}

\Address{$^\dag$~Applied and Computational Mathematics Division,
Information Technology Laboratory,\\
\hphantom{$^\dag$}~National Institute of Standards and Technology,
Gaithersburg, Maryland, USA}
\EmailD{\href{mailto:hcohl@nist.gov}{hcohl@nist.gov}}

\Address{$^\ddag$~Department of Mathematics, University of Auckland, 38 Princes Str., Auckland, New Zealand}
\URLaddressD{\url{http://www.math.auckland.ac.nz/~hcoh001/}}

\ArticleDates{Received January 19, 2011, in f\/inal form May 04, 2011;  Published online May 24, 2011}

\Abstract{For integral representations of associated Legendre functions
in terms of mo\-di\-f\/ied Bessel functions, we establish justif\/ication
for dif\/ferentiation under the integral sign with respect to parameters.
With this justif\/ication, derivatives for associated Legendre functions of
the f\/irst and second kind with respect to the degree are evaluated at
odd-half-integer degrees, for general complex-orders, and derivatives
with respect to the order are evaluated at integer-orders,
for general complex-degrees.  We also discuss the properties of the
complex function $f:\C\setminus\{-1,1\}\to\C$ given by
$f(z)=z/(\sqrt{z+1}\sqrt{z-1})$.}


\Classification{31B05; 31B10; 33B10; 33B15; 33C05; 33C10}

\section{Introduction}
\label{Introduction}

This paper is a continuation of work which is presented in %Cohl (2010)
\cite{Cohl}.  In~\cite{Cohl}, formulae were presented for derivatives of associated Legendre functions of the f\/irst
kind $P_\nu^\mu$ and the second kind $Q_\nu^\mu$ with respect to their parameters, namely
the degree $\nu$ and the order $\mu$, valid on the complex $z$-plane cut along the real axis from
$-\infty$ to $1$ (see discussion at the end of Section~\ref{ParameterderivativeformulasfromInut}).
The strategy applied in~\cite{Cohl} was to dif\/ferentiate integral representations of associated
Legendre functions, which were given in terms of modif\/ied Bessel functions of the f\/irst and
second kind, with respect to their parameters.
The derivatives of the integrands, for the integral representations of associated Legendre functions
given in~\cite{Cohl}, which include the derivatives with respect to the order evaluated at
integer-orders for modif\/ied Bessel functions of the f\/irst and second kind, are well known
(see for instance \cite[\S~3.2.3]{MOS}). % in Magnus, Oberhettinger \& Soni (1966) \cite{MOS}).

Unfortunately, in \cite{Cohl}, no
justif\/ication for dif\/ferentiation under the integral sign of the chosen integral representations
of associated Legendre functions is given.  In this paper, we give justif\/ication for dif\/ferentiation
under the integral sign for the integral representations of associated Legendre functions given
in~\cite{Cohl} and hence complete our proof for the validity of the parameter dif\/ferentiation
formulae given therein.  The parameter dif\/ferentiation formulae given in~\cite{Cohl}
are derivatives for associated Legendre functions of the f\/irst and second kind with respect
to the degree, evaluated at odd-half-integer degrees, for general complex-orders, and for
derivatives with respect to the order evaluated at integer-orders, for general complex-degrees.

There has been recent interest in the literature for tabulating closed-form
expressions of derivatives with respect to parameters for special functions
(see for instance~\cite{Brych}). %Brychkov (2008) \cite{Brych}).
Concerning derivatives with respect to parameters for associated Legendre
functions, some formulae relating to these derivatives have been previously
noted (see \cite[\S~4.4.3]{MOS}), % in Magnus, Oberhettinger \& Soni (1966) \cite{MOS}),
and there has been recent activity for solving open problems in this area
\cite{Brych, Brychkov, Cohl,Szmy0,%Szmy4,Szmy0b,
Szmy2,Szmy2b,Szmy3, Szmy1,Szmy5}.
For an extensive list of physical applications
for derivatives with respect to parameters for associated Legendre functions
with integer-order and integer-degree see %Szmytkowski (2009)
\cite{Szmy1}.
For an interesting application for derivatives of associated Legendre functions
evaluated at odd-half-integer degrees see %Szmytkowski (2011)
\cite{Szmy5}.

This paper is organized as follows.
In Section~\ref{TheWhippleformulaeforassociatedLegendrefunctions} we present a
description of a map on a subset of the complex plane which leads to the Whipple formulae
for associated Legendre functions.
In Section~\ref{Justificationfordifferentiationundertheintegralsign} we give justif\/ication for
dif\/ferentiation under the integral sign for the integral representations of associated
Legendre functions given in~\cite{Cohl}.
In Appendix~\ref{Propertiesofthefunctionzmapstofraczsqrtz1z1},
we investigate
properties of the complex function $z\mapsto {\frac{z}{\sqrt{z^2-1}}}$.

Throughout this paper we use the following conventions.  First
$ \sum\limits_{n=i}^{j}a_n=0$ for all $a_1,a_2,\ldots\in\C,$ and $i,j\in\Z$ with $j<i$.
Secondly,
for any expression
of the form $(z^2-1)^\alpha$, read this
as
\[
(z^2-1)^\alpha:=(z+1)^\alpha(z-1)^\alpha,
\]
for any f\/ixed $\alpha\in\C$ and $z\in\C\setminus\{-1,1\},$ and principal
branches are chosen
(see for instance \cite[\S~4.2(iv)]{NIST}). % of Olver {\it et al.} (2011) \cite{NIST}).

\section{The Whipple formulae for associated Legendre functions}
\label{TheWhippleformulaeforassociatedLegendrefunctions}

There is a transformation over an open subset of the complex plane which is
particularly useful in studying associated Legendre functions
(see \cite{Abra} and \cite{Hob}). %Abramowitz \& Stegun (1972) \cite{Abra} and Hobson (1955) \cite{Hob}).
This transformation, which is valid on a certain domain of the complex numbers,
accomplishes the following
\begin{gather}
\cosh z \leftrightarrow \coth w,\qquad
\coth z \leftrightarrow \cosh w,\qquad
\sinh z \leftrightarrow (\sinh w)^{-1}.
\label{transwhip}
\end{gather}
This transformation is accomplished using the map $w:{\mathfrak D}\to\C$, with
\[
{\mathfrak D}:=\C\setminus
\bigl\{z\in\C: \mathrm{Re}\, z\le 0 \mathrm{\ and\ } \mathrm{Im}\, z=2\pi n,
\ n\in\Z \bigr\},
\]
and $w$ def\/ined by
\begin{gather}
w(z):=\log\coth\frac{z}{2}.
\label{map}
\end{gather}
The map $w$ is periodic with period $2\pi i$ and is locally injective.
The map $w$ restricted to ${\mathfrak D}\cap\{z\in\C:-\pi < \mbox{Im}\, z <\pi\}$ is
verif\/ied to be an involution.
The transformation~(\ref{transwhip}) is
the restriction of the mapping~$w$
to this restricted domain.

This transformation is particularly useful for certain
associated Legendre functions such as toroidal harmonics
(see~\cite{CRTB,CT}),
%(see Cohl {\it et al.}~(2001) \cite{CRTB}, Cohl \& Tohline (1999) \cite{CT}),
associated
Legendre functions of the f\/irst and second kind with odd-half-integer degree and
integer-order, and for other associated Legendre functions which one might
encounter in potential theory.
The real argument
of toroidal harmonics
naturally occur in $(1,\infty)$, and these are the simultaneous ranges
of both the real hyperbolic cosine and cotangent functions.
One application of this map occurs with the Whipple
formulae for associated Legendre functions
\cite{CTRS,Whip}
%(Whipple (1917) \cite{Whip}, Cohl {\it et al.}~(2000) \cite{CTRS})
under index (degree and order) interchange.  See for instance  \cite[(8.2.7) and (8.2.8)]{Abra},
%(8.2.7) and (8.2.8) in Abramowitz \& Stegun (1972) \cite{Abra},
namely
\begin{subequations}
\begin{gather}
P_{-\mu-1/2}^{-\nu-1/2}
\biggl(\frac{z}{\sqrt{z^2-1}}\biggr)=
\sqrt{\frac{2}{\pi}}
\frac{(z^2-1)^{1/4}e^{-i\mu\pi}}
{\Gamma(\nu+\mu+1)}Q_\nu^\mu(z),
\label{whipple}
\end{gather}
and
\begin{gather}
Q_{-\mu-1/2}^{-\nu-1/2}
\biggl(\frac{z}{\sqrt{z^2-1}}\biggr)=
-i (\pi/2)^{1/2} \Gamma(-\nu-\mu) \left(z^2-1\right)^{1/4}e^{-i\nu\pi}
P_\nu^\mu(z),
\label{whippleb}
\end{gather}
\end{subequations}
which are valid for $\mathrm{Re}\, z>0$
and for all complex $\nu$ and $\mu,$ except where the
functions are not def\/ined.

\section{Justif\/ication for dif\/ferentiation under the integral sign}
\label{Justificationfordifferentiationundertheintegralsign}

In this section, we present and derive formulae for parameter derivatives of associated
Legendre functions of the f\/irst kind $P_\nu^\mu$ and the second
kind~$Q_\nu^\mu$, with respect to their parameters, namely the degree $\nu$ and the
order~$\mu$.
We cover parameter
derivatives of
associated Legendre functions for
argument $z\in\C\setminus(-\infty,1]$.

We incorporate derivatives with respect to order evaluated at integer-orders
for modif\/ied Bessel functions
(see \cite[\S~3.1.3, \S~3.2.3, and \S~3.3.3]{MOS}) %\S 3.1.3, \S 3.2.3, and \S 3.3.3 in Magnus, Oberhettinger \& Soni (1966) \cite{MOS})
to compute derivatives with respect to the degree and the order
of associated Legendre functions.
Below we apply these results through certain integral
representations
of associated Legendre functions
in terms of modif\/ied
Bessel functions.
Modif\/ied Bessel functions of the f\/irst and second kind respectively can be def\/ined for
all $\nu\in\C$ (see for instance \cite[\S~3.7]{Watson}) %\S3.7 in Watson (1944) \cite{Watson})
by
\[
I_\nu(z):=\sum_{m=0}^\infty \frac{(z/2)^{\nu+2m}}{m!\Gamma(\nu+m+1)},
\]
and
\[
K_\nu(z):=\frac{\pi}{2}\frac{I_{-\nu}(z)-I_\nu(z)}{\sin\pi\nu}.
\]
For $\nu=n\in\N_0:=\{0,1,2,\ldots\}$, the f\/irst def\/inition yields
\[
I_n(z)=I_{-n}(z).
\]
It may be verif\/ied that
\[
K_n(z)=\lim_{\nu\to n}K_\nu(z)
\]
is well def\/ined.
The modif\/ied Bessel function of the second kind is commonly referred to as
a~Macdonald function.

The strategy applied in this section is to use integral representations of associated
Legendre functions, expressed in terms of modif\/ied Bessel functions, and
justify dif\/ferentiation under the integral sign with respect to the relevant parameters.

\subsection[Parameter derivative formulas from $K_\nu(t)$]{Parameter derivative formulas from $\boldsymbol{K_\nu(t)}$}

It follows from \cite[(6.628.7)]{Grad}  %Gradshteyn \& Ryzhik (2007) (6.628.7) \cite{Grad}
(see also \cite[(2.16.6.3)]{Prud}) % Prudnikov {\it et al.}~(1988) (2.16.6.3) \cite{Prud})
that
\begin{gather}
  \int_0^\infty e^{-zt}K_\nu(t)t^{\mu-1/2}dt = \sqrt{\frac{\pi}{2}}
\Gamma\left(\mu-\nu+\frac12\right)\Gamma\left(\mu+\nu+\frac12\right)
\left(z^2-1\right)^{-\mu/2}
P_{\nu-1/2}^{-\mu}(z)\nonumber\\
\phantom{\int_0^\infty e^{-zt}K_\nu(t)t^{\mu-1/2}dt}{}  = \Gamma\left(\mu-\nu+\frac12\right)\left(z^2-1\right)^{-\mu/2-1/4}e^{-i\pi\nu}
Q_{\mu-1/2}^\nu\left(\frac{z}{\sqrt{z^2-1}}\right),
\label{IK}
\end{gather}
where we used the Whipple formulae (\ref{whipple}), for
$\mbox{Re}\,z>-1$ and $\mbox{Re}\,\mu>|\mbox{Re}\,\nu|-1/2$,
and the second line follows from (\ref{whippleb})
with the identity $P_\nu^\mu(z)=P_{-\nu-1}^\mu(z)$ (see \cite[(8.2.1)]{Abra}).
% in Abramowitz \& Stegun (1972) \cite{Abra}).
We would like to generate an analytical expression for the
derivative of the associated Legendre function
of the second kind with respect to its
order, evaluated at integer-orders.  In order to do this our strategy is to
solve the above integral expression for the associated
Legendre function
of the second kind,
dif\/ferentiate with respect to the order, evaluate at integer-orders, and take
advantage of the corresponding formula for dif\/ferentiation with respect to the order
for modif\/ied Bessel functions of the second kind
(see \cite[\S~3.2.3]{MOS}). %in Magnus, Oberhettinger \& Soni (1966) \cite{MOS}).
Using the expression for the associated Legendre
function of the second kind in~(\ref{IK}), we solve for $Q_{\nu-1/2}^\mu(z)$ and re-express using the map in~(\ref{map}).  This gives us the expression
\begin{gather}
Q_{\nu-1/2}^\mu(z)=\frac{\left(z^2-1\right)^{-\nu/2-1/4}e^{i\pi\mu}}{\Gamma\left(\nu-\mu+\frac12\right)}
\int_0^\infty \exp\left(\frac{-zt}{\sqrt{z^2-1}}\right)K_\mu(t)t^{\nu-1/2}dt.
\label{Qintprepwithrespecttomacdonald}
\end{gather}
By (\ref{IK}) and (\ref{whippleb}), the integral on the right-hand
side of (\ref{Qintprepwithrespecttomacdonald}) is convergent
for $\mbox{Re}\, \frac{z}{\sqrt{z+1}\sqrt{z-1}} >-1$.  By
Proposition~\ref{cotangentfunctioncomplexplane} in
Appendix~\ref{Propertiesofthefunctionzmapstofraczsqrtz1z1}, we have for all
$z\in\C\setminus[-1,1]$, this is true.  Hence the
above integral representation~(\ref{Qintprepwithrespecttomacdonald})
for the associated Legendre function of the second kind
is valid for all
$z\in\C\setminus[-1,1]$
and $\mbox{Re}\,\mu>|\mbox{Re}\,\nu|-1/2$.

In order to justify dif\/ferentiation under the
integral sign we use the following well-known corollary of the bounded convergence theorem
(cf.~\cite[\S~8.2]{Langreal}). %in Lang (1993) \cite{Langreal}).
\begin{proposition}
Let $(X,\mu)$ be a measure space, $U\subset\R$ open and $f:X\times U\to\R$ a function.  Suppose
\begin{enumerate}\itemsep=0pt
\item[$1)$] for all $y\in U$ the function $x\mapsto f(x,y)$ is measurable,
\item[$2)$] $\frac{\partial f}{\partial y}(x,y)$ exists for all $(x,y)\in X\times U$,
\item[$3)$] there exists $g\in{\mathcal L}^1(X)$
such that $\left|\frac{\partial f}{\partial y}(x,y)\right|\le g(x)$
for all $(x,y)\in X\times U$.
\end{enumerate}
Then the function $y\mapsto\int_X f(x,y) d\mu(x)$ is differentiable on $U$ and
\[
\frac{d}{dy}\left(\int_X f(x,y)d\mu(x)\right)=\int_X\frac{\partial f}{\partial y}(x,y)d\mu(x).
\]
\label{differentiateunderintegralsign}
\end{proposition}

 We call $g$ a ${\mathcal L}^1$-majorant.

We wish to dif\/ferentiate (\ref{Qintprepwithrespecttomacdonald})
with respect to the order $\mu$ and evaluate at $\mu_0=\pm m$, where $m\in\N_0$.
The derivative of the modif\/ied Bessel
function of the second kind with respect to its order
(see \cite[\S~3.2.3]{MOS}) %in Magnus, Oberhettinger \& Soni (1966) \cite{MOS})
is given by
\begin{gather}
\left[
\frac{\partial}{\partial\mu}
K_\mu(t)
\right]_{\mu=\pm m}=\pm m!
\sum_{k=0}^{m-1}
\frac{2^{m-1-k}}{k!(m-k)}t^{k-m}K_k(t)
\label{dKdn}
\end{gather}
(see for instance  \cite[(1.14.2.2)]{Brych}). %in Brychkov (2008) \cite{Brych}).
For a f\/ixed $t$, $K_\mu(t)$ is an even function of $\mu\in\R$
(see \cite[(10.27.3)]{NIST}), %in Olver {\it et al.} (2010) \cite{NIST}),
 i.e.
\[
K_{-\mu}(t)=K_{\mu}(t),
\]
and for $\mu\in[0,\infty)$, $K_\mu(t)$ is a strictly increasing function of $\mu$.
For a f\/ixed $t$, $\partial K_\mu(t)/\partial\mu$ is an odd
function of $\mu\in\R$ and for $\mu\in[0,\infty)$, and $\partial K_\mu(t)/\partial\mu$
is also a strictly increasing function of~$\mu$.  Using (\ref{dKdn}) we can make the following
estimate{\samepage
\begin{gather*}
\left|\frac{\partial}{\partial\mu}K_\mu(t)\right|<\left.
\frac{\partial K_\tau(t)}{\partial\tau}\right|_{\tau=\pm(m+1)},
%\label{dKdmuestimate}
\end{gather*}
for all $\mu\in(\mu_0-1,\mu_0+1)$.}

To justify dif\/ferentiation under the integral sign in
(\ref{Qintprepwithrespecttomacdonald}), with respect to $\mu$, evaluated at $\mu_0$,
we use Proposition
\ref{differentiateunderintegralsign}.
If we f\/ix $z$ and $\nu$, the integrand of
(\ref{Qintprepwithrespecttomacdonald}) can be given by the function
$f:\R\times(0,\infty)\to\C$ def\/ined by
\[
f(\mu,t):=\exp\left(\frac{-zt}{\sqrt{z^2-1}}\right)t^{\nu-1/2}K_\mu(t).
\]
Since $\partial K_\mu(t)/\partial\mu$ is a strictly increasing function
of $\mu\in[0,\infty),$ we have for all $\mu\in(\mu_0-1,\mu_0+1)$
\begin{gather*}
\left|\frac{\partial f}{\partial\mu}(\mu,t)\right| =
\exp\left[\mbox{Re}\left(\frac{-zt}{\sqrt{z^2-1}}\right)\right]t^{{\rm Re}\,\nu-1/2}
\left|
\frac{\partial}{\partial\mu}
K_\mu(t)
\right|\nonumber\\
\hphantom{\left|\frac{\partial f}{\partial\mu}(\mu,t)\right|}{}
 <
\exp\left[\mbox{Re}\left(\frac{-zt}{\sqrt{z^2-1}}\right)\right]t^{{\rm Re}\,\nu-1/2}
\left|\left[
\frac{\partial}{\partial\tau}
K_\tau(t)
\right]_{\tau=\pm (m+1)}\right|\nonumber\\
\hphantom{\left|\frac{\partial f}{\partial\mu}(\mu,t)\right|}{}
 =
\exp\left[\mbox{Re}\left(\frac{-zt}{\sqrt{z^2-1}}\right)\right]t^{{\rm Re}\,\nu-1/2}
\left|\left[
\frac{\partial}{\partial\tau}
K_\tau(t)
\right]_{\tau=m+1}\right| \nonumber\\
\hphantom{\left|\frac{\partial f}{\partial\mu}(\mu,t)\right|}{}
 \le
\exp\left[\mbox{Re}\left(\frac{-zt}{\sqrt{z^2-1}}\right)\right]t^{{\rm Re}\,\nu-1/2}
(m+1)!
\sum_{k=0}^{m}
\frac{2^{m-k}}{k!(m+1-k)}t^{k-m-1}K_k(t)
\nonumber\\
\hphantom{\left|\frac{\partial f}{\partial\mu}(\mu,t)\right|}{}
 \le
\exp\left[\mbox{Re}\left(\frac{-zt}{\sqrt{z^2-1}}\right)\right]t^{{\rm Re}\,\nu-1/2}
(m+1)!2^mt^{-1}K_m(t)=:g(t),\nonumber
\end{gather*}
where we used (\ref{dKdn}) and the fact that $K_k(t)\le K_m(t)$ for all
$k\in\{0,\ldots,m-1\}$.
Then $g$ is a ${\mathcal L}^1$-majorant for the derivative of the integrand, since
the integral (\ref{Qintprepwithrespecttomacdonald}) converges for
$\mbox{Re}\,(z/\sqrt{z^2-1})>-1$ and
$\mbox{Re}\,\nu>m-1/2$.

The conditions for dif\/ferentiating under the integral sign have been satisf\/ied and we can
re-write (\ref{Qintprepwithrespecttomacdonald}) as
\begin{gather*}
\left[\frac{\partial}{\partial \mu}Q_{\nu-1/2}^\mu(z)\right]_{\mu=\pm m} =
\left(z^2-1\right)^{-\nu/2-1/4}
\left[
\frac{\partial}{\partial\mu}
\frac{e^{i\pi\mu}}{\Gamma\left(\nu-\mu+\frac12\right)}
\right]_{\mu=\pm m}
%\label{derivativewithrespecttomuofQnumhlaf}
\\
\hphantom{\left[\frac{\partial}{\partial \mu}Q_{\nu-1/2}^\mu(z)\right]_{\mu=\pm m} = }{}
\times\int_0^\infty \exp\left(\frac{-zt}{\sqrt{z^2-1}}\right)K_{\pm m}(t)t^{\nu-1/2}dt
\nonumber\\
\hphantom{\left[\frac{\partial}{\partial \mu}Q_{\nu-1/2}^\mu(z)\right]_{\mu=\pm m} = }{} +
\frac{\left(z^2-1\right)^{-\nu/2-1/4}(-1)^m}{\Gamma\left(\nu\mp m+\frac12\right)}\nonumber\\
\hphantom{\left[\frac{\partial}{\partial \mu}Q_{\nu-1/2}^\mu(z)\right]_{\mu=\pm m} = }{}
\times\int_0^\infty \exp\left(\frac{-zt}{\sqrt{z^2-1}}\right)t^{\nu-1/2}
\left[
\frac{\partial}{\partial\mu}
K_\mu(t)
\right]_{\mu=\pm m}dt.
\nonumber
\end{gather*}
The derivative from the f\/irst term is given as
\[
\left[
\frac{\partial}{\partial\mu}
\frac{e^{i\pi\mu}}{\Gamma\left(\nu-\mu+\frac12\right)}
\right]_{\mu=\pm m}=
\frac{(-1)^m}{\Gamma\left(\nu\mp m+\frac12\right)}
\left[i\pi+\psi\left(\nu\mp m+\frac12\right) \right],
\]
where the $\psi$ is the digamma function
def\/ined in terms of the
derivative of the gamma function,
\[
\frac{d}{dz}\Gamma(z):=\psi(z)\Gamma(z),
\]
for $z\in\C\setminus(-\N_0)$.

Substituting these expressions for the derivatives into the two
integrals and using the map in~(\ref{map}) to re-evaluate
these integrals in terms of associated
Legendre functions gives
the following general expression for the derivative of the
associated Legendre function
of the second kind with respect to its order evaluated at
integer-orders as
\begin{gather}
\frac{\Gamma(\nu\mp m+\frac12)}{\Gamma(\nu-m+\frac12)}\left[\frac{\partial}{\partial \mu}Q_{\nu-1/2}^\mu(z)\right]_{\mu=\pm m} =
\left[i\pi+\psi\left(\nu\mp m+\frac12\right) \right] Q_{\nu-1/2}^m(z)\label{firstformula} \\
\hphantom{\frac{\Gamma(\nu\mp m+\frac12)}{\Gamma(\nu-m+\frac12)}\left[\frac{\partial}{\partial \mu}Q_{\nu-1/2}^\mu(z)\right]_{\mu=\pm m} =}{}
\pm m!\sum_{k=0}^{m-1}\frac{(-1)^{k-m}\left(z^2-1\right)^{(k-m)/2}}{k!(m-k)2^{k-m+1}}
Q_{\nu+k-m-1/2}^k(z).\nonumber
\end{gather}

If we start with the expression for the associated
Legendre function of
the f\/irst kind in (\ref{IK}) and solve for $P_{\nu-1/2}^{-\mu}(z)$ we have
\begin{gather}
P_{\nu-1/2}^{-\mu}(z)=\sqrt{\frac{2}{\pi}}
\frac{\left(z^2-1\right)^{\mu/2}}
{\Gamma\left(\mu-\nu+\frac12\right)\Gamma\left(\mu+\nu+\frac12\right)}
\int_0^\infty e^{-zt}K_\nu(t)t^{\mu-1/2}dt.
\label{PK}
\end{gather}

To justify dif\/ferentiation under the integral sign in
(\ref{PK}), with respect to $\nu$, evaluated at $\nu=\pm n$,
where $n\in\N_0$, we use as similar argument as in (\ref{Qintprepwithrespecttomacdonald})
only with modif\/ication $\mu\mapsto\nu$ and $m\mapsto n$.  The same
modif\/ied ${\mathcal L}^1$-majorant will work for the derivative of this integrand, since
the integral~(\ref{PK}) converges for
${\displaystyle \mbox{Re}\,z>-1}$ and
$\mbox{Re}\,\nu>|\mbox{Re}\,\mu|-1/2$.

The conditions for dif\/ferentiating under the integral sign have been satisf\/ied and we can
re-write (\ref{PK}) as
\begin{gather}
\left[\frac{\partial}{\partial \nu}P_{\nu-1/2}^{-\mu}
(z)\right]_{\nu=\pm n} =
\sqrt{\frac{2}{\pi}}\left(z^2-1\right)^{\mu/2}
\left[\frac{\partial}{\partial \nu}\frac{1}
{\Gamma\left(\mu-\nu+\frac12\right)\Gamma\left(\mu+\nu+\frac12\right)}
\right]_{\nu=\pm n}\label{secondderivativeofintegralforPnun}\\
\hphantom{\left[\frac{\partial}{\partial \nu}P_{\nu-1/2}^{-\mu} (z)\right]_{\nu=\pm n} =}{}
\times\int_0^\infty e^{-zt}K_{\pm n}(t)t^{\mu-1/2}dt  +
\sqrt{\frac{2}{\pi}}
\frac{\left(z^2-1\right)^{\mu/2}}
{\Gamma\left(\mu\mp n+\frac12\right)\Gamma\left(\mu\pm n+\frac12\right)}\nonumber\\
\hphantom{\left[\frac{\partial}{\partial \nu}P_{\nu-1/2}^{-\mu} (z)\right]_{\nu=\pm n} =}{}
{}\times\int_0^\infty e^{-zt}t^{\mu-1/2}
\left[\frac{\partial}{\partial \nu}K_\nu(t)\right]_{\nu=\pm n}dt.
\nonumber
\end{gather}
The derivative from the f\/irst term in (\ref{secondderivativeofintegralforPnun})
is given as
\[
\left[
\frac{\partial}{\partial\nu}
\frac{1}{\Gamma\left(\mu-\nu+\frac12\right)\Gamma\left(\mu+\nu+\frac12\right)}
\right]_{\nu=\pm n}=
\frac{\psi\left(\mu\mp n+\frac12\right)-\psi\left(\mu\pm n+\frac12\right)}
{\Gamma\left(\mu\pm n+\frac12\right)\Gamma\left(\mu\mp n+\frac12\right)}.
\]
Substituting this expression for the derivative and that given
in (\ref{dKdn}) yields the following general expression for the
derivative of the associated Legendre function
of the f\/irst kind with
respect to its degree evaluated at odd-half-integer degrees as
\begin{gather*}
\pm\left[\frac{\partial}{\partial \nu}P_{\nu-1/2}^{-\mu}
(z)\right]_{\nu=\pm n} =
\left[
\psi\left(\mu-n+\frac12\right)-\psi\left(\mu+n+\frac12\right)
\right]
P_{n-1/2}^{-\mu}(z)\nonumber\\
\qquad{} +\frac{n!}{\Gamma\left(\mu+n+\frac12\right)}
\sum_{k=0}^{n-1}
\frac{\Gamma\left(\mu-n+2k+\frac12\right)\left(z^2-1\right)^{(n-k)/2}}
{k!(n-k)2^{k-n+1}}
P_{k-1/2}^{-\mu+n-k}(z).\nonumber
\end{gather*}
If one makes a global replacement $-\mu\mapsto\mu$, using
the properties of gamma and
digamma functions, this result reduces to
\begin{gather}
\pm\left[\frac{\partial}{\partial \nu}P_{\nu-1/2}^{\mu}
(z)\right]_{\nu=\pm n} =
\left[
\psi\left(\mu+n+\frac12\right)-\psi\left(\mu-n+\frac12\right)
\right]
P_{n-1/2}^{\mu}(z)\nonumber\\
\qquad{} +n!  \Gamma\left(\mu-n+\frac12\right)
\sum_{k=0}^{n-1}
\frac{\left(z^2-1\right)^{(n-k)/2}}
{\Gamma\left(\mu+n-2k+\frac12\right)k!(n-k)2^{k-n+1}}
P_{k-1/2}^{\mu+n-k}(z).
\label{secondformula}
\end{gather}

\subsection[Parameter derivative formulas from $I_\nu(t)$]{Parameter derivative formulas from $\boldsymbol{I_\nu(t)}$}
\label{ParameterderivativeformulasfromInut}

Starting this time with \cite[(6.624.5)]{Grad} % Gradshteyn \& Ryzhik (2007) (6.624.5) \cite{Grad}
(see also \cite[(2.15.3.2)]{Prud}),
%Prudnikov {\it et al.}~(1988) (2.15.3.2) \cite{Prud}),
 we have
for $\mbox{Re}\,z>1$ and $\mbox{Re}\,\mu>-\mbox{Re}\,\nu-1/2$,
\begin{gather}
 \int_0^\infty e^{-zt}I_\nu(t)t^{\mu-1/2}dt =
\sqrt{\frac{2}{\pi}}
e^{-i\pi\mu}\left(z^2-1\right)^{-\mu/2}
Q_{\nu-1/2}^\mu(z)\nonumber\\
\hphantom{\int_0^\infty e^{-zt}I_\nu(t)t^{\mu-1/2}dt}{}
 = \Gamma\left(\mu+\nu+\frac12\right)\left(z^2-1\right)^{-\mu/2-1/4}
P_{\mu-1/2}^{-\nu}\left(\frac{z}{\sqrt{z^2-1}}\right),
\label{II}
\end{gather}
where we used again the Whipple formulae (\ref{whipple}).

We will use this particular integral representation of associated Legendre functions
to compute certain derivatives of the associated Legendre
functions with respect to the degree and the order.
We start with the integral representation of the
associated Legendre
function
of the second kind~(\ref{II}), namely
\begin{gather}
Q_{\nu-1/2}^\mu(z)
=\sqrt{\frac{\pi}{2}}e^{i\pi\mu}\left(z^2-1\right)^{\mu/2}
\int_0^\infty e^{-zt} t^{\mu-1/2}
I_\nu(t) dt.
\label{QintegralforI}
\end{gather}

To justify dif\/ferentiation under the integral sign in
(\ref{QintegralforI}), with respect to $\nu$, evaluated at $\nu_0=\pm n$,
where $n\in\N:=\{1,2,3,\ldots\}$, we use again Proposition~\ref{differentiateunderintegralsign}.
If we f\/ix $z$ and $\mu$, the integrand of
(\ref{QintegralforI})
can be given by the function
$f:\R\times(0,\infty)\to\C$ def\/ined by
\[
f(\nu,t):=e^{-zt}t^{\mu-1/2}I_\nu(t).
\]
We use the following integral representation for the derivative with
respect to the order of the modif\/ied Bessel function of the f\/irst kind
(see \cite[(75)]{ApelblatKravitsky}) %(75) in Apelblat \& Kravitsky (1985) \cite{ApelblatKravitsky})
\begin{gather}
\frac{\partial I_\nu(t)}{\partial\nu}=-\nu\int_0^t K_0(t-x)I_\nu(x)x^{-1}dx.
\label{integreetinodrvoiideselIwithnu}
\end{gather}
Let $\delta\in(0,1)$ and $M>2$.
Consider $g:(0,\infty)\to[0,\infty)$ def\/ined by
\[
g(t):=Me^{-t\Rz}t^{\Rmu-1/2}\int_0^t K_0(t-x) I_\delta(x) x^{-1}dx.
\]
Using
(\ref{integreetinodrvoiideselIwithnu})
we have
for all $\nu\in(\delta,M)$
\begin{gather*}
\left|\frac{\partial f(\nu,t)}{\partial\nu}\right| =
e^{-t\Rz }t^{\Rmu-1/2}\left|\frac{\partial I_\nu(t)}{\partial\nu}\right|
 = \nu e^{-t\Rz}t^{\Rmu-1/2}\int_0^t K_0(t-x)I_\nu(x)x^{-1}dx \\
\hphantom{\left|\frac{\partial f(\nu,t)}{\partial\nu}\right|}{}
 \le Me^{-t\Rz}t^{\Rmu-1/2}\int_0^t K_0(t-x)I_\delta(x)x^{-1}dx
 = g(t),
\end{gather*}
since for f\/ixed $t$, $\nu\mapsto I_\nu(t)$ is strictly decreasing.  Now we show that
$g\in\mathcal{L}^1$.  The integral of $g$ over its domain is
\[
\int_0^\infty g(t)dt=M\int_0^\infty e^{-t\Rz}t^{\Rmu-1/2}\int_0^t K_0(t-x) I_\delta(x)
x^{-1}dx dt.
\]
Making a change of variables in the integral, $(x,t)\mapsto(x,y)$
with $y=t-x$, yields
\[
\int_0^\infty g(t)dt=M\int_0^\infty e^{-y\Rz}K_0(y)\int_0^\infty
e^{-x\Rz}(x+y)^{\Rmu-1/2}x^{-1}I_\delta(x)dx   dy.
\]
First we show that $g$ is integrable in a neighborhood of zero.
Suppose $\mbox{Re}\,\mu-1/2<0$, $x,y\in(0,1]$ and $a\in(0,1)$. Then
\[
(x+y)^{\Rmu-1/2}=(x+y)^{-a}(x+y)^{\Rmu-1/2+a}\le y^{-a}\max\left(2^{\Rmu-1/2+a},
x^{\Rmu-1/2+a}\right).
\]
Since $K_0(y)\sim-\log(y)$ \cite[(10.30.3)]{NIST}  %((10.30.3) in Olver {\it et al.} (2010) \cite{NIST})
it follows that
\[
\int_0^1 K_0(y)y^{-a}dy<\infty.
\]
Furthermore since $I_\delta(x)\sim(x/2)^\delta/\Gamma(\delta+1)$
\cite[(10.30.1)]{NIST} %((10.30.1) in Olver {\it et al.} (2010) \cite{NIST})
it follows that
\[
\int_0^1 I_\delta(x)x^{-1} dx<\infty.
\]
Now we show that
\begin{gather}
\int_0^1 I_\delta(x)x^{\Rmu-1/2+a-1}dx,
\label{integralImuhalfa1}
\end{gather}
is convergent if $\mbox{Re}\,\mu-1/2+a+\delta>0$.
If we def\/ine
\[
 \epsilon:=\frac{\mbox{Re}\,\mu+\nu_0+\frac12}{3}>0,
\]
then $\mbox{Re}\,\mu=-\nu_0-1/2+3\epsilon$.  Therefore if we take
$a:=1-\epsilon$ and $\delta:=\nu_0-\epsilon<\nu_0$ then
\[
\mbox{Re}\,\mu-\frac12+a+\delta=\epsilon>0,
\]
and hence (\ref{integralImuhalfa1}) is convergent and thus $g$ is integrable near the origin.
If $\Rmu-1/2\ge 0$ then similarly $g$ is integrable near the origin.

Now we show that $g$ is integrable.
Suppose $\mbox{Re}\,\mu-1/2>0$. Then
\[
(x+y)^{{\rm Re}\,\mu-1/2}\le \left[2 \max(x,y))\right]^{{\rm Re}\,\mu-1/2}=2^{{\rm Re}\,\mu-1/2}\max(x^{{\rm Re}\,\mu-1/2},y^{{\rm Re}\,\mu-1/2})
\]
for all $x,y\ge 0$.
For $y\to\infty$ one has
$K_\nu(y)\sim\sqrt{\pi/(2y)}e^{-y}$ \cite[p.~250]{Olver}. % (p.~250 in Olver (1997) \cite{Olver}).
Hence it follows that
\[
\int_1^\infty K_0(y)e^{-y\Rz}y^{{\rm Re}\,\mu-1/2}dy<\infty,
\]
and
\[
\int_1^\infty K_0(y)e^{-y\Rz }dy<\infty.
\]
Furthermore since for $x\to\infty$, $I_\delta(x)\sim e^x/\sqrt{2\pi x}$
\cite[p.~251]{Olver} %(p.~251 in Olver (1997) \cite{Olver})
it follows that
\[
\int_1^\infty e^{-x\Rz} I_\delta(x)x^{{\rm Re}\,\mu-3/2} dx<\infty,
\]
and
\[
\int_1^\infty e^{-x\Rz}I_\delta(x)x^{-1} dx<\infty.
\]
If $\Rmu-1/2\le0$ then similarly $g$ is integrable.

Therefore $g$ is a
${\mathcal L}^1$-majorant for the derivative with respect to $\nu$ of the integrand
in (\ref{QintegralforI}). It is unclear whether dif\/ferentiation under the integral sign is
also possible for $\nu_0=0$.  However,
we show below that our derived results for derivatives with respect to the
degree for asso\-cia\-ted Legendre functions match up with the to be derived
results for degree $\nu=0$.  Relatively
little is known about the properties of Bessel functions in relation to
operations (dif\/ferentiation and integration) with respect to their order
(cf.~\cite{ApelblatKravitsky}).  %Apelblat \& Kravitsky (1985) \cite{ApelblatKravitsky}).

Dif\/ferentiating with respect to the degree $\nu$ and evaluating at $\nu=\pm n$, where
$n\in\N$, one obtains
\begin{gather}
\left[\frac{\partial}{\partial \nu} Q_{\nu-1/2}^\mu(z)
\right]_{\nu=\pm n}
=\sqrt{\frac{\pi}{2}}e^{i\pi\mu}\left(z^2-1\right)^{\mu/2}
\int_0^\infty e^{-zt} t^{\mu-1/2}
\left[
\frac{\partial}{\partial \nu} I_\nu(t)\right]_{\nu=\pm n}dt.
\label{dQn}
\end{gather}
The derivative of the modif\/ied Bessel function of the f\/irst kind
with respect to the order evaluated at integer-orders
(\ref{dQn}) (see \cite[\S~3.2.3]{MOS}) %\S 3.2.3 of Magnus, Oberhettinger \& Soni (1966) \cite{MOS})
is given by
\begin{gather}
\left[\frac{\partial}{\partial \nu} I_\nu(t)\right]_{\nu=\pm n}=
(-1)^{n+1}K_n(t)\pm n!\sum_{k=0}^{n-1}
\frac{(-1)^{k-n}}{k!(n-k)}
\frac{t^{k-n}}{2^{k-n+1}}
I_k(t)
\label{dIdn}
\end{gather}
(see for instance \cite[(1.13.2.1)]{Brych}). %(1.13.2.1) in Brychkov (2008) \cite{Brych}).

Inserting (\ref{dIdn}) into (\ref{dQn}) and using
(\ref{IK}) and (\ref{II}), we obtain the following general expression for the
derivative of the associated Legendre
function of the second kind with respect to its degree
evaluated at odd-half-integer degrees as
\begin{gather}
\left[
\frac{\partial}{\partial \nu} Q_{\nu-1/2}^\mu(z)
\right]_{\nu=\pm n}
=-\sqrt{\frac{\pi}{2}}e^{i\pi\mu}
\Gamma\left(\mu-n+\frac12\right)\left(z^2-1\right)^{-1/4}Q_{\mu-1/2}^n
\left(\frac{z}{\sqrt{z^2-1}}\right)\nonumber\\
\hphantom{\left[\frac{\partial}{\partial \nu} Q_{\nu-1/2}^\mu(z)\right]_{\nu=\pm n}=}{}
\pm n!
\sum_{k=0}^{n-1}\frac{\left(z^2-1\right)^{(n-k)/2}}{2^{k-n+1}k!(n-k)}
Q_{k-1/2}^{\mu+k-n}(z).
\label{derivofQwithresepecttomuevaluatpmn}
\end{gather}
Note that
\[
\left[
\frac{\partial}{\partial \nu} Q_{\nu-1/2}^\mu(z)
\right]_{\nu=0}=
-\sqrt{\frac{\pi}{2}}e^{i\pi\mu}\Gamma\left(\mu+\frac12\right)
\left(z^2-1\right)^{-1/4}
Q_{\mu-1/2}\left(\frac{z}{\sqrt{z^2-1}}\right),
\]
by \cite{MOS}. %Magnus, Oberhettinger \& Soni (1966) \cite{MOS}.
Therefore (\ref{derivofQwithresepecttomuevaluatpmn}) is also valid if $\nu=0$.

Finally, we obtain a formula for the derivative with
respect to the order for the associated
Legendre function of the f\/irst kind evaluated at integer-orders.
In order to do this we use the integral expression for the
associated Legendre function of the f\/irst kind given by (\ref{II}) and
the map given in~(\ref{map}) to convert to the appropriate
argument.  Now use the negative-order condition for associated Legendre
functions
of the f\/irst kind (see for example \cite[(22)]{CTRS}) %(22) in Cohl {\it et al.}~(2000) \cite{CTRS})
to convert to a~positive order, namely
\begin{gather}
 P_{\nu-1/2}^\mu(z) = \frac{2}{\pi}e^{-i\mu\pi}\sin(\mu\pi)Q_{\nu-1/2}^\mu(z)\nonumber
\\
\hphantom{P_{\nu-1/2}^\mu(z) =}{}
+\frac{(z^2-1)^{-\nu/2-1/4}}{\Gamma(\nu-\mu+\frac12)}\int_0^\infty
\exp\left(\frac{-zt}{\sqrt{z^2-1}} \right) I_\mu(t)t^{\nu-1/2}dt.
\label{PnuhalfmuzintegralforImu}
\end{gather}
By~(\ref{II}) and~(\ref{whippleb}), the integral on the right-hand
side of~(\ref{PnuhalfmuzintegralforImu}) converges
for $\mbox{Re}\, \frac{z}{\sqrt{z+1}\sqrt{z-1}} >1$.
By Proposition \ref{cotangentfunctioncomplexplane} in
Appendix~\ref{Propertiesofthefunctionzmapstofraczsqrtz1z1}, we only
have $\mbox{Re}\, \frac{z}{\sqrt{z+1}\sqrt{z-1}} >0$
for $z\in\C\setminus[-1,1]$. Therefore the above integral
representation~(\ref{PnuhalfmuzintegralforImu}) for the associated
Legendre function of the f\/irst kind will not be valid for the full
region $z\in\C\setminus[-1,1]$, but only for a doubly-connected open
subset of $\C$ which is symmetric about the real and imaginary axes
which includes the real segment $z\in(1,\infty)$
and has boundary in quadrant I, given by the curve
$\mbox{Re}\, \frac{z}{\sqrt{z+1}\sqrt{z-1}} =1$.  For a~detailed discussion of this curve,
see the end of Appendix~\ref{Propertiesofthefunctionzmapstofraczsqrtz1z1}.

To justify dif\/ferentiation under the integral sign in
(\ref{PnuhalfmuzintegralforImu}), with respect to $\mu$, evaluated at $\mu=\pm m$,
where $m\in\N$, we use as similar argument as in
(\ref{QintegralforI}),
only with modif\/ication $\nu\mapsto\mu$ and $n\mapsto m$.  The same
modif\/ied ${\mathcal L}^1$-majorant will work for the derivative of this integrand, since
the integral (\ref{PnuhalfmuzintegralforImu}) converges for
$\mbox{Re}\,(z/\sqrt{z^2-1})>1$ and
$\mbox{Re}\,\mu>-\mbox{Re}\,\nu-1/2$.
Since we were unable to justify dif\/ferentiation under the integral for $\nu=0$ before,
the case for dif\/ferentiation under the integral
(\ref{PnuhalfmuzintegralforImu}) with respect to $\mu$ evaluated at
$\mu=0$ remains open.  However,
below we show that our derived results for derivatives with respect to the
order for associated Legendre functions match up to previously established results
in the literature for order $\mu=0$.

Dif\/ferentiating both sides of the
resulting expression with respect to the order $\mu$ and eva\-luating
at $\mu=\pm m,$ where $m\in\N$ yields
\begin{gather*}
  \left[\frac{\partial}{\partial\mu} P_{\nu-1/2}^\mu(z) \right]_{\mu=\pm m} =
2Q_{\nu-1/2}^{\pm m}(z)+\left(z^2-1\right)^{-\nu/2-1/4}
\nonumber\\
\qquad{}
\times
\left\{
\frac{\partial}{\partial\mu}
\left[
\Gamma\left(\nu-\mu+\frac12\right)
\right]^{-1}
\right\}_{\mu=\pm m}
\int_0^\infty \exp\left(\frac{-zt}{\sqrt{z^2-1}} \right)
I_{\pm m}(t)t^{\nu-1/2}dt
\nonumber
\\
\qquad{} +\frac{\left(z^2-1\right)^{-\nu/2-1/4}}
{\Gamma\left(\nu\mp m+\frac12\right)}
\int_0^\infty \exp\left( \frac{-zt}{\sqrt{z^2-1}}\right)
t^{\nu-1/2}
\left[ \frac{\partial}{\partial\mu} I_\mu(t) \right]_{\mu=\pm m}dt.
\nonumber
\end{gather*}
The derivative of the reciprocal of the gamma function reduces to
\[
{\displaystyle \left\{
\frac{\partial}{\partial\mu}
\left[
\Gamma\left(\nu-\mu+\frac12\right)
\right]^{-1}
\right\}_{\mu=\pm m}=
\frac{\psi\left(\nu\mp m+\frac12\right)}{\Gamma\left(\nu\mp m+\frac12\right)}.}
\]
The derivative
with respect to order for the modif\/ied Bessel function of the f\/irst
kind is given in~(\ref{dIdn}).  The integrals are easily obtained by
applying the map given by~(\ref{map}) as necessary to
(\ref{IK}) and~(\ref{II}).  Hence by also using standard properties of
associated Legendre, gamma, and digamma functions we obtain the
following compact form
\begin{gather}
 \frac{\Gamma(\nu\mp m+\frac12)}{\Gamma(\nu-m+\frac12)}
\left[\frac{\partial}{\partial\mu} P_{\nu-1/2}^\mu(z) \right]_{\mu=\pm m} =
Q_{\nu-1/2}^m(z)+\psi\left(\nu\mp m+\frac12\right)
P_{\nu-1/2}^m(z)\label{differPwithresepctototomu}\\
 \hphantom{\frac{\Gamma(\nu\mp m+\frac12)}{\Gamma(\nu-m+\frac12)}
\left[\frac{\partial}{\partial\mu} P_{\nu-1/2}^\mu(z) \right]_{\mu=\pm m} = }{} \pm m!
\sum_{k=0}^{m-1}
\frac{(-1)^{k-m}\left(z^2-1\right)^{(k-m)/2}}{2^{k-m+1}k!(m-k)}
P_{\nu+k-m-1/2}^k(z).
\nonumber
\end{gather}
Note that
\[
\left[\frac{\partial}{\partial\mu} P_{\nu-1/2}^\mu(z) \right]_{\mu=0}=
Q_{\nu-1/2}(z)+\psi\left(\nu+\frac12\right)P_{\nu-1/2}(z),
\]
by \cite[\S~4.4.3]{MOS}. %of Magnus, Oberhettinger \& Soni (1966) \cite{MOS}.
So
(\ref{differPwithresepctototomu}) is also valid if~$\mu=0$.

The integral representations
(\ref{Qintprepwithrespecttomacdonald}),
(\ref{PK}),
(\ref{QintegralforI}), and
(\ref{PnuhalfmuzintegralforImu}),
which are used to obtain the parameter derivative formulae for associated Legendre functions
of the f\/irst and second kind presented in this paper (\ref{firstformula}),
(\ref{secondformula}), (\ref{derivofQwithresepecttomuevaluatpmn}), and (\ref{differPwithresepctototomu}),
are each convergent, in terms of the argument $z$, in their own
specif\/ic regions of the complex plane.
The presented parameter derivative formulae are given in terms of f\/inite sums over associated
Legendre functions which are analytic functions on the domain $\C\setminus(-\infty,1]$.
Therefore, these formulae provide an analytic continuation
for the parameter derivatives to the domain given by the cut plane with argument
$z\in\C\setminus(-\infty,1]$.

\appendix

\section[Properties of the function $z\mapsto z/\sqrt{z^2-1}$]{Properties of the function $\boldsymbol{z\mapsto z/\sqrt{z^2-1}}$}
\label{Propertiesofthefunctionzmapstofraczsqrtz1z1}

In this paper we make use of integral representations of associated Legendre functions,
namely~(\ref{IK}) and~(\ref{II}), and frequently take advantage of the
Whipple formulae~(\ref{whipple}) and~(\ref{whippleb}).  The Whipple formulae
directly relate Legendre functions of the f\/irst and second kind evaluated at
arguments~$z$ and~$z/\sqrt{z^2-1}$ respectively. Hence it is useful in conjunction
with the Whipple formulae, to understand the mapping properties
of~$z\mapsto z/\sqrt{z^2-1}$.  In particular one would like to know the
behavior of the real part of this function in regard to domains of convergence
for the integral representations of associated Legendre functions which are used.

\begin{proposition}\label{cotangentfunctioncomplexplane}
Define the function $f:\C\setminus\{-1,1\}\to\C$ by
\[
f(z)=\frac{z}{\sqrt{z^2-1}}:=\frac{z}{\sqrt{z+1}\sqrt{z-1}},
\]
where the principal branch of the square roots are chosen.
This function $f$ has the following properties.
\begin{enumerate}\itemsep=0pt
\item[$1.$] $f\bigr\vert_{\C\setminus[-1,1]}$ is even and $f\bigr\vert_{(-1,1)\pm i0}$ is odd,
where $ \pm i0:=i\lim\limits_{x\to 0^{\pm}}x$.

\item[$2.$] The sets $(0,1)\pm i0$ and $(-1,0)\pm i0$ are mapped onto
$i\left\{ \begin{matrix}
(-\infty,0)\\[1pt]
(0,\infty)
\end{matrix}\right\}$
and
$i\left\{ \begin{matrix}
 (0,\infty)\\[1pt]
 (-\infty,0)
\end{matrix} \right\}$
respectively, where
$  \pm i\infty:=i\lim\limits_{x\to\pm\infty}x$.

\item[$3.$] The sets $i(-\infty,0)$ and $i(0,\infty)$ are both mapped to $(0,1)$.

\item[$4.$] $f(0\pm i0)=0$.

\item[$5.$] If $z\in\C\setminus [-1,1]$ then ${\rm Re} \,f(z) >0$.
\end{enumerate}
\end{proposition}

\begin{proof}
When $z\ne 0$ and the exponent $w$ is any complex number, then
$z^w$ is def\/ined by the equation
\[
z^w:=\exp(w \log z),
\]
where the exponential function can be def\/ined over the entire complex plane using the power
series def\/inition
\[
\exp(z):=\sum_{n=0}^\infty \frac{z^n}{n!},
\]
and the logarithmic function is def\/ined for points $z=re^{i\arg z}$, with $r>0$, as
\[
\log z := \log r + i\arg z.
\]
Recall that if $z\in\C\setminus\{0\}$, then $\arg z$ (often referred to as the argument, amplitude or phase)
is given by the angle measured from the positive real axis to the vector representing $z$.
The angle is positive if measured anticlockwise and we choose the $\arg z\in(-\pi,\pi).$
The principal branch of the square root $\sqrt{z}$
(with branch cut along $(-\infty,0]$) is given by that unique branch of the square root
which is non-negative for $z\in(0,\infty)$.
Using this branch of the square root, the product $\sqrt{z+1}\sqrt{z-1}$ is
well-def\/ined and continuous in $z\in(-\infty,-1)$.  A branch cut along $[-1,1]$ is
chosen for $f$ which is analytic in $\C\setminus[-1,1]$.
Note that
\[
\arg (\sqrt{w})=\frac12 \arg w.
\]

If $z\in\C$ and $\mbox{Im}\,z>0$ then
\[
\arg (-(z\pm 1))=-\pi+\arg (z\pm 1),
\]
so
\[
\arg \left(\sqrt{-(z\pm 1)}\right)
=-\frac{\pi}{2}+\arg \left( \sqrt{z\pm 1}\right),
\]
and we have
\[
\sqrt{-(z\pm 1)}=-i\sqrt{z\pm 1}.
\]
Hence
\[
f(-z)=\frac{-z}{i^2\sqrt{z+1}\sqrt{z-1}}=f(z).
\]
Similarly if $\mbox{Im}\,z<0$ then
\[
\sqrt{-(z\pm 1)}=i\sqrt{z\pm 1},
\]
and we have the same result.

Let $x>1$. Then
\[
\arg \sqrt{-(x\pm 1)}=\frac{\pi}{2},
\]
so
\[
f(-x)=\frac{-x}{\sqrt{-(x-1)}\sqrt{-(x+1)}}=\frac{x}{\sqrt{x+1}\sqrt{x-1}}=f(x).
\]
Therefore $f\bigr\vert_{\C\setminus[-1,1]}$ is even.

Let $\arg   z\in(-\pi,\pi)$. For $\arg  z\gtrless 0$,
\[
f(z)=\frac{\mp i z}{\sqrt{1+z}\sqrt{1-z}},
\]
since $z-1=e^{\pm i\pi}(1-z)$. If $x\in(0,1)$, then
\[
f(x\pm i0)=\frac{\mp ix}{\sqrt{1+x}\sqrt{1-x}},
\]
and
\[
f(-x\pm i0)=\frac{\pm ix}{\sqrt{1+x}\sqrt{1-x}}=-f(x\pm i0).
\]
Moreover, $f(0\pm i0)=0$.
Therefore $f\bigr\vert_{(-1,1)\pm i0}$ maps to the imaginary axis and is an odd
function of $x$.

If $x\in(0,\infty)$ then
\[
f(ix)=\frac{ix}{\sqrt{ix+1}\sqrt{ix-1}}=\frac{x}{\sqrt{1+x^2}},
\]
and
\[
f(-ix)=\frac{-ix}{\sqrt{-ix+1}\sqrt{-ix-1}}=\frac{x}{\sqrt{1+x^2}},
\]
so $f$ maps both the positive and negative imaginary axes to the real
interval $(0,1)$. Clearly $f(0)=0$.  This completes the proof of 1, 2, 3 and 4.

\begin{figure}[t]
\centering
\includegraphics[width=70mm,angle=90]{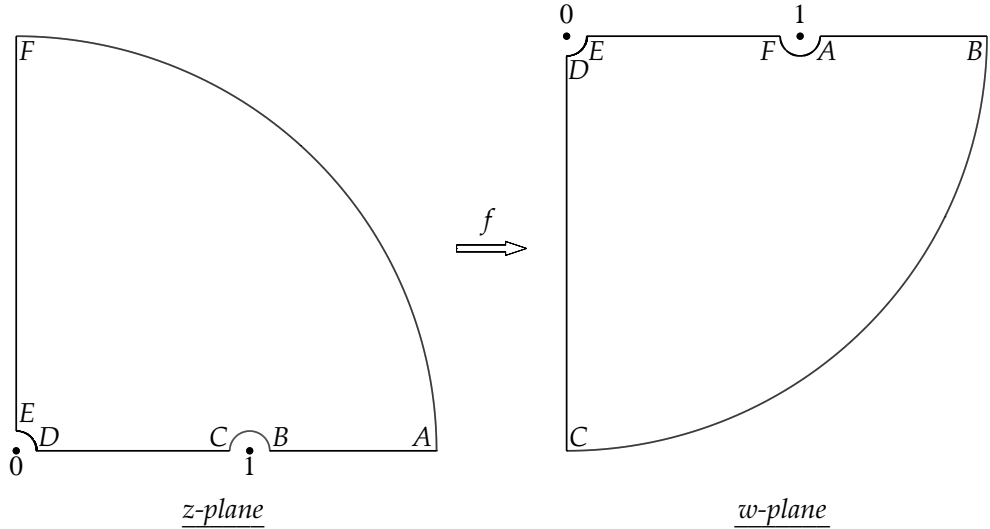}

\caption{This f\/igure shows how the function $f:\C\setminus\{-1,1\}\to\C$ def\/ined by
$f(z):=z/(\sqrt{z+1}\sqrt{z-1})$ conformally maps quadrant I into quadrant~IV.}
\label{Fig:conformal}
\end{figure}

Before we prove 5 we f\/irst show that $f$ maps quadrant I into quadrant IV.  The
derivative $f^\prime(z)=-(z+1)^{-3/2}(z-1)^{-3/2}$ is nowhere zero, and therefore
$w=f(z)$ represents a conformal map of $\C\setminus\{-1,1\}$
(see for instance \cite[\S~1.9(iv)]{NIST}). %\S 1.9(iv) of Olver {\it et al.} (2010) \cite{NIST}).
Consider the closed
contour represented in the $z$-plane of Fig.~\ref{Fig:conformal}.  In order to
study the mapping properties of quadrant I by the map $f$, we look at the
behavior of the map within and on a closed contour in
quadrant~I through a limiting process:
Take the radii of the semi-circular segment $BC$ and the quarter-circular segment
$DE$ tending towards zero, the radius of the quarter-circular segment $FA$ tending
towards inf\/inity, and match these segments to the straight line
segments $AB$, $CD$, and $EF$ continuously.

The straight line segments are treated f\/irst, followed by the treatment of the circular segments.
For $AB$, $z\in(1,\infty)$ and
therefore $w(z)=1/\sqrt{1-(1/z^2)}\in(1,\infty)$. Therefore as $z\to A$, $w(z)\to 1{+}$ and as
$z\to B$, $w(z)\to +\infty$.
On $CD$, $z\in(0,1)$ so
$w(z)=-i/\sqrt{(1/z^2)-1}\in -i(0,\infty)$.
So as $z\to C$, $w(z)\to -i\infty$ and as
$z\to D$, $w(z)\to -i0$.
On $EF$, $z\in i(0,\infty)$ thus
$w(z)=z/\sqrt{1+z^2}\in(0,1)$ and as $z\to E$, $w(z)\to 0{+}$ and as
$z\to F$, $w(z)\to 1{-}.$
For the semi-circular segment $BC$, $z$ is near $1$ and we write $z=1+\zeta$.
Consider $\zeta=\epsilon e^{i\phi}$
with $\phi\in[0,\pi]$.
Through the binomial expansion we can see that
\begin{gather}
w(\zeta)=\frac{1}{\sqrt{2\zeta}}+O\big(\sqrt{\zeta}\big)=\frac{1}{\sqrt{2\epsilon}}e^{-i\phi/2}+O\big(\sqrt{\epsilon}\big).
\label{asymptoticsclosetoone}
\end{gather}
Therefore as $z\to B$, $w(z)\to+\infty$ and as
$z\to C$, $w(z)\to-i\infty.$
For the quarter-circular segment $DE$, $z$ is near $0$ and we write $z=\zeta=\epsilon e^{i\phi}$
with $\phi\in[0,\pi/2]$.
Through the binomial expansion we have
\[
w(\zeta)=-i\zeta+O\big(\zeta^2\big)=\epsilon e^{i(\phi-\pi/2)}+O\big(\epsilon^2\big),
\]
and therefore as $z\to D$, $w(z)\to-i0$ and as $z\to E$, $w(z)\to 0{+}$.
For the quarter-circular segment $FA$, we write $z=\zeta$
and consider $\zeta=R e^{i\phi}$
with $\phi\in[0,\pi/2]$ with $R$ chosen suf\/f\/iciently large.
Through the binomial expansion we have{\samepage
\begin{gather}
w(\zeta)=1+\frac{1}{2\zeta^2}+O\big(\zeta^{-4}\big)
=1+\frac{1}{2R^2}e^{-2i\phi}+O\big(R^{-4}\big),
\label{asymptoticsinfinity}
\end{gather}
and therefore as $z\to F$, $w(z)\to 1{-}$.}

Hence the closed contour and its interior region in the $z$-plane
are conformally mapped into the closed contour and its interior region
in the $w$-plane shown in Fig.~\ref{Fig:conformal}.
Through by the limiting process described above, we see that $f$ maps quadrant I into quadrant IV.

Due to the evenness of $f$, quadrants I \& III are
mapped to quadrants IV, and quadrants II \& IV are mapped to quadrant I. Therefore
if $z\in\C\setminus [-1,1]$ then $\mbox{\rm Re}\,f(z)>0 $.
This completes the proof of 5.
\end{proof}

\noindent
{\bf Note.} An anonymous referee has suggested an alternate proof that
$\mbox{Re}\, f(z)>0$ for all $\C\setminus[-1,1]$ using the minimum modulus
principle.

\begin{proof}
Consider the function $g:\C\setminus\{-1,1\}\to\C$ def\/ined by
$g(z):=\exp f(z).$  A simple computation gives the modulus $|g(z)|=\exp\mbox{Re}\,f(z)$.
By~(\ref{asymptoticsinfinity}), $|g(z)|\to e$ as $|z|\to\infty$.
Let $x\in(-1,1)$ and
consider $z=x+i\epsilon$ as $\epsilon\to 0^\pm$, then
\[
f(z)=\frac{\mp ix}{\sqrt{1+x}\sqrt{1-x}}+O(\epsilon),
\]
so $|g(z)|\to 1$ in $(-1,1)$ (clearly $g$ is non-constant).   In a small neighborhood $E$ of $z=1$,
and by property~1, in a small neighborhood $E^\prime$ of $z=-1$,
we have through~(\ref{asymptoticsclosetoone}) that $|g(z)|> 1$ in~$E$ and~$E^\prime$.
The minimum modulus principle states (see for instance \cite[p.~147]{Silverman}) %p.147 in Silverman (1974) \cite{Silverman})
that {\em if $f$ is analytic, non-constant, and
non-vanishing in an open connected subset $G$ of $\C$ then $|f(z)|$ cannot have
a~minimum in $G$.} Let $G:=\C\setminus[-1,1]$, an open connected subset of $\C$.
Since $f$ (and hence $g$) is analytic in $G$, and the exponential function is
non-vanishing in $\C$, from the minimum modulus principle we have that
$|g(z)|>1$ in $\C\setminus[-1,1]$ and therefore $\mbox{Re}\, f(z)>0$ in
$\C\setminus[-1,1]$.  This completes the proof.
\end{proof}

The range of $f$ is
$\{z\in\C:\mbox{Re}\, z\ge 0 \mbox{\ and\ }z\ne 1\}.$
Every complex
number in the range of the function is taken twice except for elements in $(0,1)$ and on
the imaginary axis. These complex numbers are taken only once.

\begin{figure}[t]
\centering
\includegraphics[width=90mm,angle=90]{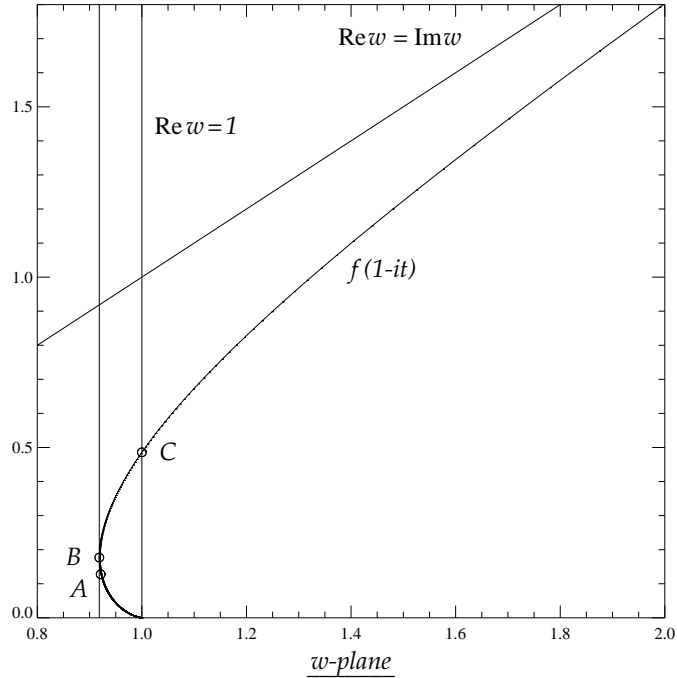}

\caption{This f\/igure shows the curve in the complex $w$-plane given
by $f(1-it)$ for $t\in(0,\infty)$.  There are three lines indicated.
One represents the line $\mbox{Re}\,w=\mbox{Im}\,w$.  The other two
represent the vertical lines passing through the points~$B$ and~$C$
respectively.}
\label{Fig:conformal2}
\end{figure}

Consider the curve $f(z)$ when $\mbox{Re}\,z=1$.
In order to illustrate the behavior of this curve, take $z=1-it$, where $t\in(0,\infty)$.
By the above discussion, we know this line segment is mapped conformally from
quadrant IV to quadrant I.
This smooth curve asymptotically approaches
the line $\mbox{Re}\, w=\mbox{Im}\, w$ from below and approaches the singularity
at unity from the left side (i.e.\ from $\mbox{Re}\,f(z)<1)$. This can be seen through
the asymptotics
(cf.~(\ref{asymptoticsclosetoone}) and~(\ref{asymptoticsinfinity})),
namely $\phi(\epsilon) = \pi-2\sqrt{2\epsilon}$
for the angle of approach to the singularity at unity as
the distance $\epsilon\to 0$ and
$\phi(R)=\pi/4-11\sqrt{2}/(64R^2)$
(where terms up to fourth order have been included)
as the radius~$R$ tends towards inf\/inity.  If $f(1-it)=x(t)+iy(t)$, then using
elementary trigonometry one can show
\begin{gather*}
x(t)=
\frac{1}{2\sqrt{t(4+t^2)}}
\left[
(1+t)\sqrt{\sqrt{4+t^2}+2}
+(t-1)\sqrt{\sqrt{4+t^2}-2}
\ \right],\nonumber\\
y(t)=
\frac{1}{2\sqrt{t(4+t^2)}}
\left[
(1-t)\sqrt{\sqrt{4+t^2}+2}
+(t+1)\sqrt{\sqrt{4+t^2}-2}
\ \right],\nonumber
\end{gather*}
and therefore
\begin{gather*}
x^2(t)=
\frac{1}{2t(4+t^2)}
\left[
(1+t^2)\sqrt{4+t^2}+3t+t^3
\right],
\nonumber\\
y^2(t)=
\frac{1}{2t(4+t^2)}
\left[
(1+t^2)\sqrt{4+t^2}-3t-t^3
\right].\nonumber
\end{gather*}
The distance squared from the origin is given by $(1+t^2)/(t\sqrt{4+t^2})$,
whose minimum occurs at $t=\sqrt{2}$ for
\[
(x,y)(\sqrt{2})=\left(
\frac12\sqrt{\frac{3\sqrt{3}+5}{3}},
\frac12\sqrt{\frac{3\sqrt{3}-5}{3}}
\right),
\]
portrayed by the point $A$ in Fig.~\ref{Fig:conformal2}.
Due to the asymptotics of the curve near the singularity at unity,
there exists a point at which the real part of this curve reaches a minimum value.
By f\/inding the minimum of $x^2(t)$,
this point can be easily obtained and is
shown in Fig.~\ref{Fig:conformal2} by point~$B$ and is given
at $t=2/\sqrt{3}$ by the point
\[
(x,y)\left(\frac{2}{\sqrt{3}}\right)=\left(
\frac34\sqrt{\frac{3}{2}},
\frac{1}{4\sqrt{2}}
\right).
\]
Similarly one can f\/ind the point $C$ in Fig.~\ref{Fig:conformal2} where
$\mbox{Re}\,f(z)=1$ to be given at $t=\sqrt{\sqrt{5}-2}$ at $y=t$.

\subsection*{Acknowledgements}

I would like to thank Dr.\ A.F.M.~ter Elst for extremely valuable discussions and
acknowledge funding for time to write this paper from the
Dean of the Faculty of Science at the University of Auckland in the form of a three
month stipend to enhance University of Auckland 2012 PBRF Performance.
I would like to express my gratitude to the anonymous referees whose helpful
comments improved this paper. I would also like to thank F.W.J.~Olver  for helpful discussions.
Part of this work was conducted while the author was a National Research Council
Research Postdoctoral Associate in the Information Technology Laboratory of the
National Institute of Standards and Technology.

\pdfbookmark[1]{References}{ref}
\LastPageEnding

\end{document}